# Relativistic constitutive modeling of inelastic deformation of continua moving in space-time

Eun-Ho Lee[1,2*]

1. Department of Mechanical Engineering, Sungkyunkwan University, Suwon, Gyeonggi-do, 16419, Republic of Korea
2. Department of Smart Fab. Technology, Sungkyunkwan University, Suwon, Gyeonggi-do, 16419, Republic of Korea

 * Corresponding author: Eun-Ho Lee, e-mail: e.h.lee@skku.edu

## Abstract

Since Einstein introduced the theory of relativity, many scientific observations have proven that it thoroughly approximates the motion of materials in space-time. Thus, some efforts have been made to expand classical continuum mechanics to relativistic continuum models to consider the relativistic effect. However, most models only consider the elastic range without accounting for inelastic deformation. A relativistic inelasticity model should provide objective measures of inelastic deformation for different observers with nonzero relative velocities over a wide speed range, even when the moving speed is close to the speed of light. This paper presents a modeling structure of the relativistic constitutive equations of the inelastic deformation of a material moving with a wide speed range. In this model, the deformation tensors of classical mechanics are expanded to four-dimensional tensors (called relativistic Cauchy–Green deformation tensors in this paper) in space-time, based on which constitutive equations of inelastic deformation are introduced. The four-dimensional tensors are objective about the homogeneous Lorentz transformation in the Minkowski space-time, and the material dissipation determined using the proposed modeling of inelastic deformation satisfies the second law of thermodynamics. In addition, an example illustrates how to apply this theory with explaining the physical meanings of the modeling. Finally, it is also shown that the proposed relativistic inelasticity modeling is collapsed to a classical inelasticity model when the speed of motion is much slower than the speed of light.

Keywords: Relativistic continuum mechanics; Thermodynamics; Constitutive equation; Inelastic deformation; Space-time

# 1. Introduction

The invariant theory is important for establishing a mechanics theory to obtain objective measures for several observers with different velocities. The Galilean invariance group is effective and provides a decent framework for this purpose when the magnitude of the relative velocity between an observer and an observed material point is much smaller than the speed of light. This assumption is called a non-relativistic condition in this paper. The modern theory of continuum mechanics has led to the development of balance laws and constitutive equations based on the Galilean invariance [1-5]. Because the Galilean group considers that time and space are independent, non-relativistic models have been developed in the Euclidean space.

In non-relativistic continuum mechanics, a smooth mapping function is defined with respect to time $t$ to describe the deformation of a body $\mathcal{B}$, which is expressed as

$$\mathbf{x} = \mathbf{x}(\mathbf{X}, t), \tag{1}$$

where $\mathbf{X}$ is the position vector ($\mathbf{X} = X^K \mathbf{e}_K$, $K = 1-3$) of a material point in the reference configuration $\kappa_R$, and changes to the current position $\mathbf{x}$ ($\mathbf{x} = x^k \mathbf{e}_k$, $k = 1-3$) in the current configuration $\kappa$ at each time $t$. $\mathbf{e}_i$ is the basis of the position vectors and $x^i$ is the contravariant component. $x^i$ is independent of $t$ in the non-relativistic mechanics. The Euclidean metric is

$$|\mathbf{x}|^2 = \mathbf{x} \cdot \mathbf{x} = \mathbf{x}^T \mathbf{I} \mathbf{x}, \tag{2}$$

where "·" denotes a vector dot product and **I** is the identity matrix. Moreover, reverse mapping can be obtained by

$$\mathbf{X} = \mathbf{X}(\mathbf{x}, t); \tag{3}$$

Thus, deformation gradient tensor **F** is expressed as

$$\mathbf{F} = \frac{\partial \mathbf{x}}{\partial \mathbf{X}}. \tag{4}$$

Because **F** and the Euclidean metric in equation (2) are objective about the Galilean transformation, objective constitutive equations can be defined under the non-relativistic conditions based on equations (1–4). This can be achieved when the magnitude of the velocity of the observed material point is much smaller than the speed of light $c$ in vacuum. When a non-relativistic assumption cannot be applied, the Galilean group cannot provide invariant measures, particularly as the magnitude of the relative velocity approaches the speed of light. Equations (1-4) should be expanded to be invariant even when the moving speed is close to the speed of light.

Lorentz transformation (Lorentz, 1899; 1904) [6-7] has provided a decent framework for invariant measures, even when the magnitude of the relative velocity is close to the speed of light [8-9]. In 1905, Einstein [10] discussed the relationship between space and time to deal with the motion of a body with a speed similar to that of light; this is known as the theory of special relativity. Briefly, he reconstructed Galilean transformations of nonrelativistic mechanics into Lorentz transformations and provided a physical meaning to the Lorentz transformations; time and space are not absolute measures and influence each other. The theory of special relativity was further refined mathematically using the concept of the Minkowski space-time [11]. The Minkowski space provides a mathematical tool to easily understand special relativity by combining the three-dimensional space and time into a four-dimensional space-time. Note that the Minkowski space is not a four-dimensional Euclidean space; it shows that time and space are a unified four-dimensional space-time continuum, which cannot be separately defined. Modern scientific observations have shown that the theory of relativity can approximate the realities of space and time.

Based on these scientific advances, efforts have been made to expand non-relativistic continuum mechanics to a relativistic continuum model by defining an invariant measured about the Lorentz transformation in space-time. The position vector in the deformed configuration in equation (1) should be redefined as a four-dimensional vector, and the Euclidean metric in equation (2) needs to be replaced by another metric to satisfy a Lorentz transformation (Møller, 1952 [12]; Toupin, 1957 [13]). Grot and Eringe [14] extended the Toupin's study [13] and introduced an elastic model that satisfies the law of thermodynamics with special relativity. In addition, the conservation of mass must be reformulated because mass changes [10,15] as the moving speed approaches the speed of light. Moreover, the motion of particles requires defining their relationship to the world velocity [16]. Once a continuum model is extended to be invariant about a Lorentz transformation, its combination with electromagnetism becomes more sophisticated mathematically [17,18]. Recent studies have considered viscoelasticity [19] and elasticity modeling in general relativity [20]. However, the aforementioned models are only in the elastic deformation range.

The inelasticity theory has been developed over the past century through experimental studies and mathematical modeling. Depending on an observer, Lagrangian [21-23] or Eulerian [24-27] modeling can be used; however, both perspectives have been established based on the Galilean transformations [28-30]. Inelasticity models need to consider the natural characteristics of inelastic deformation to define the rate of inelastic deformation. In addition, because inelastic deformation is an irreversible process that results in material dissipation, the constitutive equations should satisfy the second law of thermodynamics. Consequently, a relativistic inelasticity model should satisfy the objectivity about Lorentz transformations while considering the characteristics of inelasticity.

This study is aimed at establishing a modeling structure of the relativistic constitutive equation of the inelastic deformation of continua moving over a wide speed range, even when the speed is close to that of light. In this model, the extended four-dimensional deformation tensors defined by Grot and Eringe [14], which are objective about homogeneous Lorentz transformations, are used. The extended deformation tensors are called relativistic Cauchy–Green deformation tensors in this study, because their forms and meanings can be considered as extensions of the Cauchy–green tensors to which the relativistic effect is added. Thus, the constitutive equations of the inelastic deformation part are formulated based on the above-mentioned four-dimensional tensors. Inelastic deformation is based on extending the traditional decomposition of the deformation gradient, inspired by Lee's theorem [21]. Therefore, the proposed model has a contact point with non-relativistic inelasticity models. This study shows how to define the direction of the relativistic inelastic deformation rate, which satisfies the second law of thermodynamics. A scalar function controlling the magnitude of the relativistic inelastic deformation is also defined to satisfy the yield consistency, hardening, and rate of the stress tensor. Consequently, the proposed model has the advantage of being able to use previously introduced inelasticity theories after modifications to consider the relativity effect. After the mathematical derivation is completed, the physical meanings of the formulations are explained using a simple two-dimensional Minkowski space example by comparing two observers in different frames. This example also illustrates how to apply the proposed theory, and shows that even though the two observers have different observations, they have invariant measures for inelastic deformation. The relationship between the proposed and non-relativistic models is also illustrated by showing that the proposed inelasticity model can be identical to the nonrelativistic model when the speed of motion is much smaller than that of light.

The remainder of this paper is organized as follows. The kinematics of relativistic continuum mechanics are introduced inspired by [13-15] in section 2. Section 3 proposes a relativistic constitutive model of inelastic deformation. Section 4 presents an example with discussing the physical meaning of the model and the relationship between the proposed and nonrelativistic models. Finally, conclusions are presented in Section 5.

## 2. Kinematics of relativistic continuum mechanics

### 2.1. Space-time and four-vector

In the Minkowski space, time is considered as the fourth coordinate. In this frame, a four-vector $\hat{\pmb{x}}$ is defined to represent the coordinates of an 'event' of a material particle in the space-time, which is expressed as

$$\hat{\mathbf{x}} = x^\alpha \hat{\mathbf{e}}_\alpha \ (a = 1-4), \ \text{where} \ x^4 = ct \ . \tag{5}$$

Note that the superscript "^" denotes a four dimensional tensor in this paper, and the 'event' means the spatial and temporal position of a material particle in the Minkowski space. $\hat{\mathbf{x}}$ replaces the position vector $\mathbf{x}$, which is in equation (1). Reverse function $\mathbf{X}(\ )$ can be defined, similar to equation (3),

$$\mathbf{X} = \mathbf{X}(\hat{\mathbf{x}}). \tag{6}$$

$\mathbf{X}$ shares the same material particle with $\hat{\mathbf{x}}$, and physically means the three dimensional position vector of the reference configuration at rest in a frame. The Minkowski space-time defines a new metric $\boldsymbol{\eta}$ for four-vectors as

$$|\hat{\mathbf{x}}|^2 = \hat{\mathbf{x}} \cdot \hat{\mathbf{x}} = \hat{\mathbf{x}}^T \boldsymbol{\eta} \hat{\mathbf{x}} \text{ where } \boldsymbol{\eta} = \boldsymbol{\eta}^T = \begin{bmatrix} 1 & 0 & 0 & 0 \\ 0 & 1 & 0 & 0 \\ 0 & 0 & 1 & 0 \\ 0 & 0 & 0 & -1 \end{bmatrix}. \tag{7}$$

Because the Galilean transformation does not keep the invariance with the Minkowski metric, the group of homogeneous Lorentz transformation $\boldsymbol{\Lambda}$ should be applied to the four-vector $\hat{\mathbf{x}}$ as

$$\hat{\mathbf{x}}^* = \boldsymbol{\Lambda}\hat{\mathbf{x}}, \tag{8a}$$

$$\boldsymbol{\Lambda}\boldsymbol{\Lambda}'^T = \boldsymbol{\Lambda}'^T\boldsymbol{\Lambda} = \mathbf{I} \text{ where } \boldsymbol{\Lambda}' = \boldsymbol{\eta}\boldsymbol{\Lambda}\boldsymbol{\eta} \tag{8b}$$

$$\boldsymbol{\Lambda}^T \boldsymbol{\eta} \boldsymbol{\Lambda} = \boldsymbol{\eta} \tag{8c}$$

$$\boldsymbol{\eta}^* = \boldsymbol{\Lambda}\boldsymbol{\eta}\boldsymbol{\Lambda}^T = \boldsymbol{\eta}, \tag{8d}$$

$$\det(\boldsymbol{\Lambda}) = \pm 1. \tag{8e}$$

$\hat{\mathbf{x}}^*$ denotes the transformed four-vector by $\boldsymbol{\Lambda}$. Based on equations (7) and (8), it can be shown that $|\hat{\mathbf{x}}|^2$ are invariant with the group of homogeneous Lorentz transformation,

$$|\hat{\mathbf{x}}|^2 = \hat{\mathbf{x}} \cdot \hat{\mathbf{x}} = \boldsymbol{\Lambda}'\hat{\mathbf{x}}^* \cdot \boldsymbol{\Lambda}'\hat{\mathbf{x}}^* = \hat{\mathbf{x}}^{*T}\boldsymbol{\eta}\hat{\mathbf{x}}^* = |\hat{\mathbf{x}}^*|^2. \tag{9}$$

The continuous sequence of the event of a particle, related to the history of a particle, draws a path curve $W$, called the word line, as shown in figure 1. Figure 1 conceptually illustrates the four-dimensional Minkowski space consisting of a flat three dimensional space ($\hat{\mathbf{e}}_1$, $\hat{\mathbf{e}}_2$, and $\hat{\mathbf{e}}_3$) and time ($\hat{\mathbf{e}}_4$). The tangent to the world line $W$ defines a unit four-vector called world velocity $\hat{\mathbf{u}}$ as

$$\hat{\mathbf{u}} = \hat{\mathbf{u}}(\hat{\mathbf{x}}) = \hat{\mathbf{u}}(\mathbf{X}(\hat{\mathbf{x}})) = \begin{bmatrix} \frac{\mathbf{v}}{c\sqrt{1-\beta^2}} \\ \frac{1}{\sqrt{1-\beta^2}} \end{bmatrix}, \text{ where } \beta = \frac{|\mathbf{v}|}{c}, \tag{10a}$$

$$\mathbf{v} = \frac{\partial \mathbf{x}}{\partial t} \text{ and } |\mathbf{v}| = \sqrt{\mathbf{v} \cdot \mathbf{v}}, \tag{10b}$$

$$|\hat{\mathbf{u}}|^2 = \hat{\mathbf{u}} \cdot \hat{\mathbf{u}} = \hat{\mathbf{u}}^T \boldsymbol{\eta} \hat{\mathbf{u}} = -1, \tag{10c}$$

$$\left[\frac{\partial \mathbf{x}}{\partial \hat{\mathbf{x}}}\right] \hat{\mathbf{u}} = \mathbf{0}. \tag{10d}$$

$\mathbf{v}$ is the particle velocity. Equation (10d) can be generalized to an invariant derivative of a function $D(\hat{\mathbf{A}})$ in special relativistic mechanics,

$$D(\hat{\mathbf{A}}) = \left[\frac{\partial \hat{\mathbf{A}}}{\partial \hat{\mathbf{x}}}\right] \hat{\mathbf{u}} = \widehat{gard}(\hat{\mathbf{A}}) \hat{\mathbf{u}}, \tag{11}$$

where $\widehat{gard}(\hat{\mathbf{A}})$ provides the partial derivatives of $\hat{\mathbf{A}}$ in the four-dimensional space-time.

Next, a typical four-vector $\hat{\mathbf{f}}$ is classified as a time-like vector when $|\hat{\mathbf{f}}|^2 < 0$, space-like vector when $|\hat{\mathbf{f}}|^2 > 0$, or null if $|\hat{\mathbf{f}}|^2 = 0$. $\hat{\mathbf{f}}$ can be generally decomposed to a time-like component $\hat{\mathbf{f}}_t$ and space like component $\hat{\mathbf{f}}_s$ [14,16], which expressed as

$$\hat{\mathbf{f}} = \hat{\mathbf{f}}_t + \hat{\mathbf{f}}_s, \tag{12a}$$

$$\hat{\mathbf{f}}_t = f_t \hat{\mathbf{u}} \text{ and } \hat{\mathbf{f}}_s \cdot \hat{\mathbf{f}}_t = 0. \tag{12b}$$

$f_t$ is a scalar value. Equation (12b) means that $\hat{\mathbf{f}}_t$ is parallel to $\hat{\mathbf{u}}$ while $\hat{\mathbf{f}}_s$ is perpendicular to $\hat{\mathbf{u}}$. The two components can be obtained by a projector tensor $\mathbf{S}$

$$\hat{\mathbf{f}}_s = \mathbf{S}\hat{\mathbf{f}} \text{ and } \hat{\mathbf{f}}_t = (\mathbf{I} - \mathbf{S})\hat{\mathbf{f}}, \tag{13a}$$

$$\mathbf{S} = \mathbf{I} + \hat{\mathbf{u}} \otimes (\boldsymbol{\eta} \hat{\mathbf{u}}), \tag{13b}$$

$$\mathbf{S}\hat{\mathbf{u}} = \hat{\mathbf{u}}^T \mathbf{S}^T = 0, \mathbf{SS} = \mathbf{S}, \text{ and } \mathbf{S}^* = \boldsymbol{\Lambda} \mathbf{S} \boldsymbol{\Lambda}'^T. \tag{13c}$$

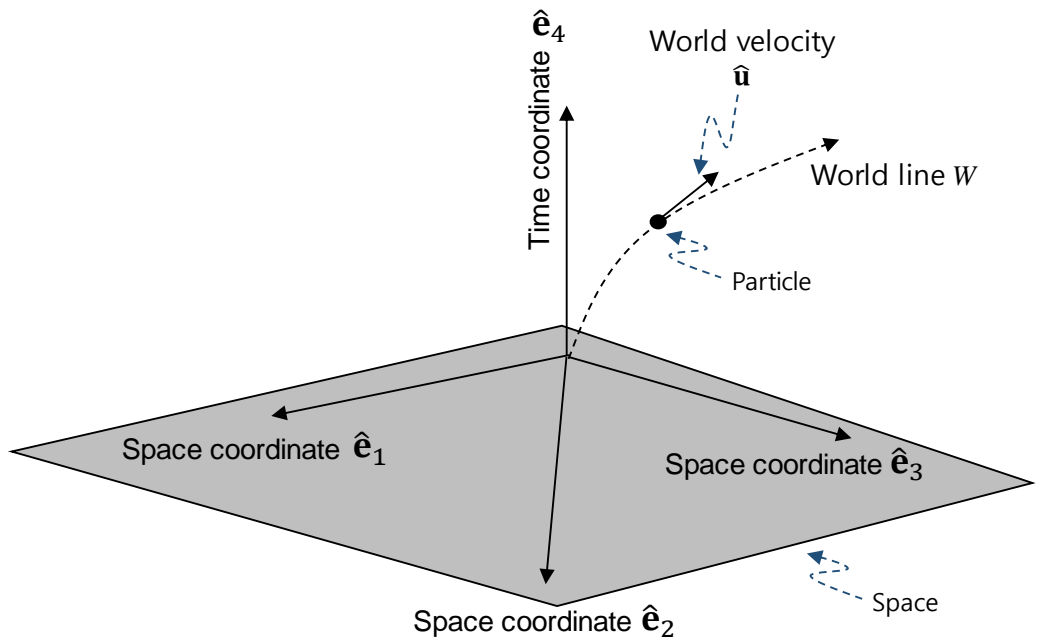

Figure 1. Concept of the relativistic kinematics in the Minkowski space

## 2.2. Deformation

Because the deformation gradient **F** of non-relativistic mechanics in equation (4) has no objectivity in the Minkowski space, a four dimensional deformation gradient tensor is defined as

$$d\hat{\mathbf{x}} = \hat{\mathbf{F}}d\mathbf{X} \text{ where } \hat{\mathbf{F}} = \frac{\partial \hat{\mathbf{x}}}{\partial \mathbf{X}}. \tag{14}$$

Based on equations (8a-8e), it can be shown that

$$d\hat{\mathbf{x}}^* = \hat{\mathbf{F}}^* d\mathbf{X} \text{ where } \hat{\mathbf{F}}^* = \mathbf{\Lambda}\hat{\mathbf{F}}. \tag{15}$$

To define deformation measures, it is useful to extract the space-like components from $d\hat{\mathbf{x}}$ as

$$d\hat{\mathbf{x}}_s = \mathbf{S}d\hat{\mathbf{x}} = \hat{\mathbf{F}}_s d\mathbf{X}, \text{ where } \hat{\mathbf{F}}_s = \mathbf{S}\hat{\mathbf{F}}, \tag{16a}$$

$$d\hat{\mathbf{x}}_s^* = \hat{\mathbf{F}}_s^* d\mathbf{X}, \text{ where } \hat{\mathbf{F}}_s^* = \mathbf{\Lambda}\hat{\mathbf{F}}_s. \tag{16b}$$

The Minkowski metric then provides

$$|d\hat{\mathbf{x}}_s|^2 = d\hat{\mathbf{x}}_s \cdot d\hat{\mathbf{x}}_s = d\hat{\mathbf{x}}_s^T \mathbf{\eta} d\hat{\mathbf{x}}_s = d\mathbf{X}^T \hat{\mathbf{C}} d\mathbf{X}, \tag{17a}$$

$$\hat{\mathbf{C}} = \hat{\mathbf{F}}_s^T \mathbf{\eta} \hat{\mathbf{F}}_s. \tag{17b}$$

They are invariant in the Mikowski space, which is shown by

$$|d\hat{\mathbf{x}}_s^*|^2 = d\hat{\mathbf{x}}_s^* \cdot d\hat{\mathbf{x}}_s^* = d\hat{\mathbf{x}}_s^{*T}\boldsymbol{\eta}d\hat{\mathbf{x}}_s^* = d\mathbf{X}^T\hat{\mathbf{C}}d\mathbf{X} = |d\hat{\mathbf{x}}_s|^2. \tag{18}$$

The reverse of $\hat{\mathbf{F}}$ also gives that

$$|d\mathbf{X}|^2 = d\mathbf{X} \cdot d\mathbf{X} = d\mathbf{X}^T d\mathbf{X} = d\hat{\mathbf{x}}_s^T \hat{\mathbf{B}}^{-1} d\hat{\mathbf{x}}_s, \tag{19a}$$

$$\hat{\mathbf{B}} = \hat{\mathbf{F}}_s \hat{\mathbf{F}}_s^T, \tag{19b}$$

$$d\mathbf{X}^T d\mathbf{X} = d\hat{\mathbf{x}}_s^{*T} \hat{\mathbf{B}}^{*-1} d\hat{\mathbf{x}}_s^*, \text{ where } \hat{\mathbf{B}}^* = \boldsymbol{\Lambda}\hat{\mathbf{F}}_s \hat{\mathbf{F}}_s^T \boldsymbol{\Lambda}^T = \boldsymbol{\Lambda}\hat{\mathbf{B}}\boldsymbol{\Lambda}^T. \tag{19c}$$

Inspired by Grot and Eringe [14], $\hat{\mathbf{C}}$ and $\hat{\mathbf{B}}$ in equations (17b) and (19b) can be considered as the relativistic right and left Cauchy-Green deformation tensors, respectively. Equations (18) and (19c) show that $\hat{\mathbf{C}}$ is not affected by the Lorentz transformation and $\hat{\mathbf{B}}$ is objective under the Lorentz transformation. Consequently, $\hat{\mathbf{C}}$ and $\hat{\mathbf{B}}$ can be used as objective deformation measures in relativistic continuum mechanics. $\hat{\mathbf{B}}$ presents information set of deformation observed by an observer, in a space-time frame $S$, who has a relative velocity for the observed material point. On the other hand, $\hat{\mathbf{C}}$ gives a deformation measure for another observer who shares the same frame (local instantaneous rest frame $S^{'}$) with the observed material point, neglecting the effect of relative velocity. $\hat{\mathbf{C}}$ and $\hat{\mathbf{B}}$ have similar forms and roles to the classical right and left Cauchy-Green deformation tensors in non-relativistic mechanics, respectively; detailed roles of $\hat{\mathbf{C}}$ and $\hat{\mathbf{B}}$ are presented in section 4.2. The two tensors $\hat{\mathbf{C}}$ and $\hat{\mathbf{B}}$ have invariants, given by

$$I_1 = \hat{\mathbf{C}}:\mathbf{I} = \hat{\mathbf{B}}:\boldsymbol{\eta}\mathbf{I}, \tag{20a}$$

$$I_2 = (\hat{\mathbf{C}}\mathbf{I}\hat{\mathbf{C}}):\mathbf{I} = (\hat{\mathbf{B}}\boldsymbol{\eta}\hat{\mathbf{B}}):\boldsymbol{\eta}\mathbf{I}. \tag{20b}$$

In this study, the operator ":" for two tensors (**A** and **B**) provides a scalar value calculated by $tr(\mathbf{AB}^T)$.

Next, it is convenient to define below four-dimensional tensor,

$$\hat{\mathbf{L}} = \widehat{gard}(\hat{\mathbf{u}}) = \frac{\partial \hat{\mathbf{u}}}{\partial \hat{\mathbf{x}}}. \tag{21}$$

Equations (11) and (16) provide

$$\mathbf{S}D(\hat{\mathbf{F}}_s) = \hat{\mathbf{L}}\hat{\mathbf{F}}_s, \tag{22}$$

and equations (17) and (22) show that

$$D(\hat{\mathbf{C}}) = 2\hat{\mathbf{F}}_s^T \hat{\mathbf{D}}_s^{\boldsymbol{\eta}} \hat{\mathbf{F}}_s, \tag{23a}$$

$$\widehat{\mathbf{D}}_s^{\eta} = \tfrac{1}{2}\left(\widehat{\mathbf{L}}_s^{\eta} + \widehat{\mathbf{L}}_s^{\eta T}\right) = \boldsymbol{\eta}\widehat{\mathbf{D}}_s\boldsymbol{\eta}^T, \text{ where } \widehat{\mathbf{D}}_s = \tfrac{1}{2}\left(\widehat{\mathbf{L}}_s + \widehat{\mathbf{L}}_s^T\right), \tag{23b}$$

$$\widehat{\mathbf{L}}_s^{\eta} = \boldsymbol{\eta}\widehat{\mathbf{L}}_s, \text{ where } \widehat{\mathbf{L}}_s = \widehat{\mathbf{L}}\mathbf{S}. \tag{23c}$$

Equations (23a-23c) show that $\widehat{\mathbf{L}}_s^{\eta}$ and $\widehat{\mathbf{D}}_s^{\eta}$ work as the relativistic velocity gradient and rate of deformation, respectively. The relativistic spin rate $\widehat{\mathbf{W}}_s^{\eta}$ is then defined by

$$\widehat{\mathbf{W}}_s^{\eta} = \widehat{\mathbf{L}}_s^{\eta} - \widehat{\mathbf{D}}_s^{\eta} = \tfrac{1}{2}\left(\widehat{\mathbf{L}}_s^{\eta} - \widehat{\mathbf{L}}_s^{\eta T}\right) = \boldsymbol{\eta}\widehat{\mathbf{W}}_s\boldsymbol{\eta}^T, \text{ where} \tag{24a}$$

$$\widehat{\mathbf{W}}_s = \widehat{\mathbf{L}}_s - \widehat{\mathbf{D}}_s = \tfrac{1}{2}\left(\widehat{\mathbf{L}}_s - \widehat{\mathbf{L}}_s^T\right). \tag{24b}$$

Volume change can be measured by the ratio between the deformed volume $dv$ in an instantaneous rest frame and the un-deformed volume $dV$. A Jacobian $j$ is defined by

$$j = \frac{dv}{dV} = \det(\begin{bmatrix}\widehat{\mathbf{F}}_s & \widehat{\mathbf{u}}\end{bmatrix}), \tag{25}$$

and the invariant derivative of $j$ is given by

$$D(j) = j \cdot tr(\widehat{\mathbf{L}}_s). \tag{26}$$

Next, the conservation of mass is replaced by conservation of the number of particles because mass is not conserved in the theory of relativity [10,14]. The conservation of particle number is based on the assumption of non-creation or indestrucibility of particles. The conservation of the particle number is given by

$$D(m_0) + m_0 \widehat{div}(\widehat{\mathbf{u}}) = 0. \tag{27}$$

$m_0$ denotes the rest frame particle number. $\widehat{div}(\ )$ denotes divergence of a four-tensor in the Minkowski space.

## 3. Relativistic inelastic deformation

### 3.1. Decomposition

Lee [21] proposed a scheme that total deformation gradient can be multiplicatively separated into elastic and inelastic parts, which has been widely applied in many works [5,31,32]. This paper introduces decompositions of the four dimensional deformation gradient, which is expressed as

$$\hat{\mathbf{F}}_s = \hat{\mathbf{F}}_s^e \mathbf{F}^p, \tag{28a}$$

$$\hat{\mathbf{F}}_s^e = \mathbf{S}\hat{\mathbf{F}}^e, \tag{28b}$$

$$\hat{\mathbf{F}}^e = \frac{\partial \hat{\mathbf{x}}}{\partial \mathbf{X}_0} \text{ and } \mathbf{F}^p = \frac{\partial \mathbf{X}_0}{\partial \mathbf{X}}. \tag{28c}$$

$\hat{\mathbf{F}}_s^e$ and $\mathbf{F}^p$ denote the elastic and inelastic parts of the deformation gradient, respectively, and $\mathbf{X}_0$ is a position vector of an intermediate configuration in a rest frame. As discussed in many studies on inelasticity modeling [33-37], defining an intermediate configuration does not have a unique way and depends on constitutive equation. Inspired by this reason, intermediate configuration vector can be defined as a four-vector in a different way, expressed as

$$\hat{\mathbf{F}}_s = \hat{\mathbf{F}}_s^{e\dagger} \boldsymbol{\eta} \hat{\mathbf{F}}^p, \tag{29a}$$

$$\hat{\mathbf{F}}_s^{e\dagger} = \mathbf{S}\hat{\mathbf{F}}^{e\dagger}, \tag{29b}$$

$$\hat{\mathbf{F}}^{e\dagger} = \frac{\partial \hat{\mathbf{x}}}{\partial \hat{\mathbf{x}}_0} \text{ and } \hat{\mathbf{F}}^p = \frac{\partial \hat{\mathbf{x}}_0}{\partial \mathbf{X}}. \tag{29c}$$

$\hat{\mathbf{x}}_0$ is a four-vector in another intermediate configuration in another frame having a relative velocity to the observed material point. Equations (28) and (29) have a relation between the two definitions as

$$\hat{\mathbf{F}}_s = \hat{\mathbf{F}}_s^{e\dagger} \boldsymbol{\eta} \hat{\mathbf{F}}^p = \hat{\mathbf{F}}_s^e \mathbf{F}^p, \tag{30a}$$

$$\hat{\mathbf{F}}^p = \hat{\mathbf{F}}_M \mathbf{F}^p, \text{ where } \hat{\mathbf{F}}_M = \frac{\partial \hat{\mathbf{x}}_0}{\partial \mathbf{X}_0}, \tag{30b}$$

$$\hat{\mathbf{F}}_s^{e\dagger} = \hat{\mathbf{F}}_s^e \hat{\mathbf{F}}_M^T, \tag{30c}$$

$$\hat{\mathbf{F}}_M^T \boldsymbol{\eta} \hat{\mathbf{F}}_M = \mathbf{I}. \tag{30d}$$

Equation (30d) means that $\hat{\mathbf{F}}_M$ does not result in additional deformation.

Next, as discussed in inelasticity theories, decomposition of $\hat{\mathbf{F}}_s$ allows that $\hat{\mathbf{L}}$ in equation (21) can be decomposed into elastic and inelastic parts

$$\hat{\mathbf{L}} = \hat{\mathbf{L}}^e + \hat{\mathbf{L}}^p, \tag{31a}$$

$$[\hat{\mathbf{L}}_s^{\boldsymbol{\eta}}]^e = \boldsymbol{\eta}\hat{\mathbf{L}}^e\mathbf{S} = \boldsymbol{\eta}\hat{\mathbf{L}}_s^e, \text{ where } \hat{\mathbf{L}}_s^e = \hat{\mathbf{L}}^e\mathbf{S}, \tag{31b}$$

$$[\hat{\mathbf{L}}_s^{\boldsymbol{\eta}}]^p = \boldsymbol{\eta}\hat{\mathbf{L}}^p\mathbf{S} = \boldsymbol{\eta}\hat{\mathbf{L}}_s^p, \text{ where } \hat{\mathbf{L}}_s^p = \hat{\mathbf{L}}^p\mathbf{S}. \tag{31c}$$

Decomposition of $\hat{\mathbf{L}}^e$ and $\hat{\mathbf{L}}^p$ requires to define a constitutive equation of inelasticity. The details of the constitutive equation of relativistic inelastic deformation are defined in section 3.3. Based on equations (23), (24), (25) and (31), $\widehat{\mathbf{D}}_s^{\eta}$, $\widehat{\mathbf{W}}_s^{\eta}$ and $j$ also can be decomposed as

$$\widehat{\mathbf{D}}_s^{\eta} = [\widehat{\mathbf{D}}_s^{\eta}]^e + [\widehat{\mathbf{D}}_s^{\eta}]^p, \text{ and } \widehat{\mathbf{W}}_s^{\eta} = [\widehat{\mathbf{W}}_s^{\eta}]^e + [\widehat{\mathbf{W}}_s^{\eta}]^p, \tag{32a}$$

$$[\widehat{\mathbf{D}}_s^{\eta}]^e = \frac{1}{2}\left([\widehat{\mathbf{L}}_s^{\eta}]^e + [\widehat{\mathbf{L}}_s^{\eta}]^{e^T}\right) = \boldsymbol{\eta}\widehat{\mathbf{D}}_s^e\boldsymbol{\eta}^T, \text{ and } [\widehat{\mathbf{D}}_s^{\eta}]^p = \frac{1}{2}\left([\widehat{\mathbf{L}}_s^{\eta}]^p + [\widehat{\mathbf{L}}_s^{\eta}]^{p^T}\right) = \boldsymbol{\eta}\widehat{\mathbf{D}}_s^p\boldsymbol{\eta}^T, \tag{32b}$$

$$\widehat{\mathbf{D}}_s^e = \frac{1}{2}\left(\widehat{\mathbf{L}}_s^e + \widehat{\mathbf{L}}_s^{e^T}\right), \text{ and } \widehat{\mathbf{D}}_s^p = \frac{1}{2}\left(\widehat{\mathbf{L}}_s^p + \widehat{\mathbf{L}}_s^{p^T}\right), \tag{32c}$$

$$[\widehat{\mathbf{W}}_s^{\eta}]^e = \frac{1}{2}\left([\widehat{\mathbf{L}}_s^{\eta}]^e - [\widehat{\mathbf{L}}_s^{\eta}]^{e^T}\right) = \boldsymbol{\eta}\widehat{\mathbf{W}}_s^e\boldsymbol{\eta}^T, \text{ and } [\widehat{\mathbf{W}}_s^{\eta}]^p = \frac{1}{2}\left([\widehat{\mathbf{L}}_s^{\eta}]^p - [\widehat{\mathbf{L}}_s^{\eta}]^{p^T}\right) = \boldsymbol{\eta}\widehat{\mathbf{W}}_s^p\boldsymbol{\eta}^T \tag{32d}$$

$$\widehat{\mathbf{W}}_s^e = \frac{1}{2}\left(\widehat{\mathbf{L}}_s^e - \widehat{\mathbf{L}}_s^{e^T}\right), \text{ and } \widehat{\mathbf{W}}_s^p = \frac{1}{2}\left(\widehat{\mathbf{L}}_s^p - \widehat{\mathbf{L}}_s^{p^T}\right), \tag{32e}$$

$$j = j^e \cdot j^p = j^e, \text{ where } j^e = \det\left(\begin{bmatrix}\widehat{\mathbf{F}}_s^e & \widehat{\mathbf{u}}\end{bmatrix}\right) \tag{32f}$$

Note that this work employs the zero dilatancy of inelastic deformation in the rest frame; $j^p = 1$. The elastic parts of $\widehat{\mathbf{C}}$ and $\widehat{\mathbf{B}}$ are given by

$$\widehat{\mathbf{C}}^e = \widehat{\mathbf{F}}_s^{e^T}\boldsymbol{\eta}\widehat{\mathbf{F}}_s^e = \left(\widehat{\mathbf{F}}_s^{e\dagger}\boldsymbol{\eta}\mathbf{M}\right)^T\left(\widehat{\mathbf{F}}_s^{e\dagger}\boldsymbol{\eta}\mathbf{M}\right), \tag{33a}$$

$$\widehat{\mathbf{B}}^e = \widehat{\mathbf{F}}_s^e\widehat{\mathbf{F}}_s^{e^T} = \widehat{\mathbf{F}}_s^{e\dagger}\boldsymbol{\eta}\widehat{\mathbf{F}}_s^{e\dagger^T}. \tag{33b}$$

The inelastic part of $\widehat{\mathbf{C}}$ and $\widehat{\mathbf{B}}$ also can be obtained by the same way using $\mathbf{F}^p$.

### 3.2. Inequality condition of material dissipation

The free energy density per particle number is defined by

$$\psi_{00} = e - h_{00}\theta, \tag{34}$$

where $e$ is the internal energy density per particle number, $h_{00}$ represents the entropy density per particle number, and $\theta$ denotes the absolute temperature. Inspired by Grot and Eringen [14], projection of the balance of momentum onto the world velocity leads to the first law of thermodynamics,

$$\widehat{\mathbf{u}} \cdot \widehat{div}(\widehat{\mathbf{T}}) = \widehat{\mathbf{f}} \cdot \widehat{\mathbf{u}}. \tag{35}$$

$\widehat{\mathbf{T}}$ denotes the energy-momentum tensor in space-time while $\widehat{\mathbf{f}}$ includes body force and energy supply density. The specific form of the energy momentum tensor $\widehat{\mathbf{T}}$ varies depending on the definition of the problem. This study focuses on inelastic deformation without considering electromagnetic interaction and heat flow in the, resulting in,

$$\widehat{\mathbf{T}} = w\hat{\mathbf{u}} \otimes \hat{\mathbf{u}} - \hat{\mathbf{t}}_s. \tag{36}$$

$w$ means energy density and $\hat{\mathbf{t}}_s$ denotes relativistic stress tensor. In this case, the components of $w$ and $\hat{\mathbf{t}}$ can be obtained by

$$w = m_0 e, \tag{37a}$$

$$\hat{\mathbf{t}}_s = -\mathbf{S}\widehat{\mathbf{T}}\mathbf{S}^{\mathrm{T}}, \tag{37b}$$

$$\hat{\mathbf{t}}_s = \hat{\mathbf{t}}_s^{\mathrm{T}}, \tag{37c}$$

$$\widehat{\mathbf{T}} = \widehat{\mathbf{T}}^{\mathrm{T}}. \tag{37d}$$

$\hat{\mathbf{f}}$ is

$$\hat{\mathbf{f}} = \begin{bmatrix} \mathbf{f} \\ h + \mathbf{f} \cdot \mathbf{v} \end{bmatrix}, \tag{38}$$

where $\mathbf{f}$ is the body force vector and $h$ denotes the internal heat generation. Substituting equation (36) into equation (35), with integration by parts, leads to that

$$\widehat{div}(w\hat{\mathbf{u}}) = \hat{\mathbf{t}}_s : \eta \widehat{grad}(\hat{\mathbf{u}}) - \hat{\mathbf{f}} \cdot \hat{\mathbf{u}}. \tag{39}$$

Next, the Clausius–Duhem inequality condition in space-time can be given by (Grot and Eringe, 1966a; Lianis, 1973)

$$\widehat{div}(\hat{\mathbf{h}}) + r \geq 0. \tag{40a}$$

$$\hat{\mathbf{h}} = h_0 \hat{\mathbf{u}} + \frac{\hat{\mathbf{q}}}{\theta} \text{ where } h_0 = m_0 h_{00}, \tag{40b}$$

$$r = \frac{r_0}{\theta} \text{ where } r_0 = \hat{\mathbf{f}} \cdot \hat{\mathbf{u}}. \tag{40c}$$

$\hat{\mathbf{h}}$ is the entropy four-vector including a time-like component $h_0 \hat{\mathbf{u}}$ and a space-like component $\frac{\hat{\mathbf{q}}}{\theta}$. $r$ represents the entropy supply. As discussed with energy momentum tensor, the heat flow is not considered under an assumption of isothermal condition in this study; $\hat{\mathbf{q}} = \mathbf{0}$ in equation (40b). Substituting equations (37a), (39) and (40b) into equation (40a) presents the inequality condition of material dissipation rate $\xi$ as below:

$$\xi = m_0[D(h_{00})\theta - D(e)] + \hat{\mathbf{t}}_s : \eta \widehat{grad}(\hat{\mathbf{u}}) \geq 0. \tag{41}$$

The invariant derivative of the free energy in equation (34) is given by

$$D(\psi_{00}) = D(e) - D(h_{00})\theta - h_{00}D(\theta). \tag{42}$$

Substituting equation (42) into (41) re-arranges the rate of material dissipation ξ

$$\xi = \hat{\mathbf{t}}_s : \boldsymbol{\eta} \widehat{grad}(\hat{\mathbf{u}}) - m_0[D(\psi_{00}) + h_{00}D(\theta)] \geq 0 \ . \tag{43}$$

In this study, $\psi_{00}$ can be given by a function of $\hat{\mathbf{C}}^e$ and $\theta$,

$$\psi_{00} = \psi_{00}(\hat{\mathbf{C}}^e, \theta). \tag{44}$$

The invariant derivative of $\psi_{00}$ is expressed as

$$D(\psi_{00}) = \frac{\partial \psi_{00}}{\partial \hat{\mathbf{C}}^e} : D(\hat{\mathbf{C}}^e) + \frac{\partial \psi_{00}}{\partial \theta} D(\theta). \tag{45}$$

With equations (43) and (45), ξ becomes

$$\xi = \left[\hat{\mathbf{t}}_s - 2m_0 \hat{\mathbf{F}}_s^e \frac{\partial \psi_{00}}{\partial \hat{\mathbf{C}}^e} \hat{\mathbf{F}}_s^{eT}\right] : \left[\hat{\mathbf{D}}_s^{\boldsymbol{\eta}}\right]^e + \hat{\mathbf{t}} : \left[\hat{\mathbf{D}}_s^{\boldsymbol{\eta}}\right]^p - m_0 \left[\frac{\partial \psi_{00}}{\partial \theta} D(\theta) + h_{00}D(\theta)\right] \geq 0 \ . \tag{46}$$

Then ξ is finalized for the relativistic inelastic deformation by

$$\xi = \hat{\mathbf{t}}_s : \left[\hat{\mathbf{D}}_s^{\boldsymbol{\eta}}\right]^p \geq 0 \ , \text{ where} \tag{47a}$$

$$\hat{\mathbf{t}}_s = 2m_0 \hat{\mathbf{F}}_s^e \frac{\partial \psi_{00}}{\partial \hat{\mathbf{C}}^e} \hat{\mathbf{F}}_s^{eT}, \tag{47b}$$

$$-\frac{\partial \psi_{00}}{\partial \theta} = h_{00} \ . \tag{47c}$$

The proposed constitutive equation for the rate of relativistic inelastic deformation $\left[\hat{\mathbf{D}}_s^{\boldsymbol{\eta}}\right]^p$ should satisfy the inequality condition of material dissipation in equation (47a).

### 3.3. Constitutive equation of relativistic inelastic deformation

Based on equation (47a), the partial derivative of ξ results in

$$\left[\hat{\mathbf{D}}_s^{\boldsymbol{\eta}}\right]^p = \frac{\partial \xi}{\partial \hat{\mathbf{t}}_s} \ . \tag{48}$$

Inspired by the flow rule [38-40], $\left[\hat{\mathbf{D}}_s^{\boldsymbol{\eta}}\right]^p$ also can be re-defined by

$$\left[\hat{\mathbf{D}}_s^{\boldsymbol{\eta}}\right]^p = \overline{\left[\hat{\mathbf{D}}_s^{\boldsymbol{\eta}}\right]}^p D(\Gamma_p), \tag{49}$$

where $\overline{\left[\hat{\mathbf{D}}_s^{\boldsymbol{\eta}}\right]}^p$ physically means the direction of $\left[\hat{\mathbf{D}}_s^{\boldsymbol{\eta}}\right]^p$ and $D(\Gamma_p)$ controls the scalar magnitude of $\left[\hat{\mathbf{D}}_s^{\boldsymbol{\eta}}\right]^p$. This equation also provides the direction of $\left[\hat{\mathbf{L}}_s^{\boldsymbol{\eta}}\right]^p$ as

$$\left[\hat{\mathbf{L}}_s^{\boldsymbol{\eta}}\right]^p = \overline{\left[\hat{\mathbf{L}}_s^{\boldsymbol{\eta}}\right]}^p D(\Gamma_p), \text{ where} \tag{50a}$$

$$\overline{[\hat{\mathbf{D}}_s^\eta]}^p = \tfrac{1}{2}(\overline{[\hat{\mathbf{L}}_s^\eta]}^p + \overline{[\hat{\mathbf{L}}_s^\eta]}^{pT}). \tag{50b}$$

The rate of material dissipation $\xi$ consist of an objective inelastic potential function $Q_p(\hat{\mathbf{t}}_s)$ and $D(\varGamma_p)$, expressed as

$$\xi = Q_p(\hat{\mathbf{t}}_s)D(\varGamma_p) = Q_p(\hat{\mathbf{t}}_s^*)D(\varGamma_p), \tag{51a}$$

$$\hat{\mathbf{t}}_s^* = \mathbf{\Lambda}\hat{\mathbf{t}}_s\mathbf{\Lambda}^T. \tag{51b}$$

Combining equations (48), (49) and (51) provides a rate of relativistic inelastic deformation, given by

$$[\hat{\mathbf{D}}_s^\eta]^p = \overline{[\hat{\mathbf{D}}_s^\eta]}^p D(\varGamma_p), \tag{52a}$$

$$\overline{[\hat{\mathbf{D}}_s^\eta]}^p = \frac{\partial Q_p}{\partial \hat{\mathbf{t}}_s}. \tag{52b}$$

Equation (52b) means that $\overline{[\hat{\mathbf{D}}_s]}^p$ is the outward normal to the surface of $Q_p(\hat{\mathbf{t}}_s)$. With equations (8) and (32b), it can be shown that $[\hat{\mathbf{D}}_s^\eta]^p$ is an objective tensor as

$$[\hat{\mathbf{D}}_s^\eta]^{p*} = \mathbf{\eta}^* \hat{\mathbf{D}}_s^{p*} \mathbf{\eta}^{T*} = \mathbf{\Lambda}'[\hat{\mathbf{D}}_s^\eta]^p \mathbf{\Lambda}'^T, \tag{53a}$$

$$\hat{\mathbf{D}}_s^{p*} = \mathbf{\Lambda}\hat{\mathbf{D}}_s^p\mathbf{\Lambda}^T. \tag{53b}$$

By setting $Q_p(\hat{\mathbf{t}}_s)$ as a convex surface, $[\hat{\mathbf{D}}_s^\eta]^p$ satisfies the inequality condition of equation (47a).

Next, motivated by consistency modeling [5, 37, 39], loading condition $F$ is defined by

$$F(\hat{\mathbf{t}}_s, \varGamma_p) = f(\hat{\mathbf{t}}_s) - t_y(\varGamma_p) = 0, \text{ where} \tag{54a}$$

$$t_y(\varGamma_p) = m_0 t_{y0}(\varGamma_p). \tag{54b}$$

$f(\hat{\mathbf{t}}_s)$ is the yield surface and $t_{y0}(\varGamma_p)$ denotes the flow stress per particle number with respect to the scalar hardening variable $\varGamma_p$. When $F(\hat{\mathbf{t}}_s, \varGamma_p) < 0$, the material undergoes pure elastic deformation condition. During inelastic deformation occurs ($F = 0$), the consistency condition satisfies that

$$D(F) = D(f) - D(t_y) = 0. \tag{55}$$

With equations (27) and (54), the consistency condition of equation (55) becomes

$$\frac{\partial f}{\partial \hat{\mathbf{t}}_s}:\hat{\mathbf{t}}_s^\nabla - m_0 \frac{\partial t_{y0}}{\partial \varGamma_p}D(\varGamma_p) + \widehat{div}(\hat{\mathbf{u}})t_y = 0. \tag{56}$$

$\hat{\mathbf{t}}_s^\nabla$ is an objective time rate of $\hat{\mathbf{t}}_s$, which can be derived based on equations (11), (23a), (27), (47b) and (51b), as

$$\hat{\mathbf{t}}_s^{\nabla} = D(\hat{\mathbf{t}}_s) + \mathbf{\Lambda}'\widehat{grad}(\mathbf{\Lambda})\hat{\mathbf{u}}\hat{\mathbf{t}}_s + \hat{\mathbf{t}}_s\widehat{grad}(\mathbf{\Lambda}^T)\hat{\mathbf{u}}\mathbf{\Lambda}'^T, \text{ where} \tag{57a}$$

$$D(\hat{\mathbf{t}}_s) = \hat{\mathbf{L}}_s^e\hat{\mathbf{t}}_s + \hat{\mathbf{t}}_s\hat{\mathbf{L}}_s^{eT} - \widehat{div}(\hat{\mathbf{u}})\hat{\mathbf{t}}_s + 4m_0\hat{\mathbf{F}}_s^e[\frac{\partial\psi_{00}}{\partial\hat{\mathbf{C}}^e\otimes\partial\hat{\mathbf{C}}^e}:\hat{\mathbf{F}}_s^{eT}(\hat{\mathbf{D}}_s^{\eta} - [\hat{\mathbf{D}}_s^{\eta}]^p)\hat{\mathbf{F}}_s^e]\hat{\mathbf{F}}_s^{eT} . \tag{57b}$$

Combining equations (56) and (57a-c) provides the value of $D(\Gamma_p)$ as

$$D(\Gamma_p) = \frac{g_1}{g_2}, \tag{58a}$$

$$g_1 = \frac{\partial f}{\partial \hat{\mathbf{t}}_s}:\{\hat{\mathbf{L}}_s^e\hat{\mathbf{t}}_s + \hat{\mathbf{t}}_s\hat{\mathbf{L}}_s^{eT} - \widehat{div}(\hat{\mathbf{u}})\hat{\mathbf{t}}_s + 4m_0\hat{\mathbf{F}}_s^e(\frac{\partial\psi_{00}}{\partial\hat{\mathbf{C}}^e\otimes\partial\hat{\mathbf{C}}^e}:\hat{\mathbf{F}}_s^{eT}\hat{\mathbf{D}}_s^{\eta}\hat{\mathbf{F}}_s^e)\hat{\mathbf{F}}_s^{eT} + \mathbf{\Lambda}'\widehat{grad}(\mathbf{\Lambda})\hat{\mathbf{u}}\hat{\mathbf{t}}_s +$$
$$\hat{\mathbf{t}}_s\widehat{grad}(\mathbf{\Lambda}^T)\hat{\mathbf{u}}\mathbf{\Lambda}'^T\} + \widehat{div}(\hat{\mathbf{u}})t_y, \tag{58b}$$

$$g_2 = m_0\frac{\partial t_{y0}}{\partial \Gamma_p} + 4m_0\frac{\partial f}{\partial \hat{\mathbf{t}}_s}:\hat{\mathbf{F}}_s^e(\frac{\partial\psi_{00}}{\partial\hat{\mathbf{C}}^e\otimes\partial\hat{\mathbf{C}}^e}:\hat{\mathbf{F}}_s^{eT}[\overline{\hat{\mathbf{D}}_s^{\eta}}]^p\hat{\mathbf{F}}_s^e)\hat{\mathbf{F}}_s^{eT} . \tag{58c}$$

Equations (52) and (58) complete the relativistic constitutive equation of inelastic deformation rate.

## 4. Example and discussions

### 4.1. Example and physical interpretation

This section discusses the physical meaning of the proposed model using the two-dimensional Minkowski space consisting of one-dimensional space axis $\hat{\mathbf{e}}_1$ and one-dimensional time axis $\hat{\mathbf{e}}_2$ for intuitive understanding. For applying the constitutive equation to specific problems, the detailed forms of yield surface $f$, hardening stress $t_y$, free energy density $\psi_{00}$, and inelastic potential function $Q_p$ should be specified. The specific forms of $f$, $t_y$, $\psi_{00}$ and $Q_p$ can vary depending on the definition of problem. Many models have been reported to discuss the detailed forms of $f$, $t_y$, $\psi_{00}$ and $Q_p$. A common form of $\psi_{00}$ is given by a function of the three invariants $I_1^e$, $I_2^e$, and $I_3^e$, expressed as

$$\psi_{00} = \psi_{00}(I_1^e, I_2^e, I_3^e), \tag{59a}$$

$$I_1^e = \hat{\mathbf{C}}^e : \mathbf{I}, \tag{59b}$$

$$I_2^e = \frac{1}{2}[I_1^{e^2} - (\hat{\mathbf{C}}^e\mathbf{I}\hat{\mathbf{C}}^e):\mathbf{I}], \tag{59c}$$

$$I_3^e = \det(\hat{\mathbf{C}}^e). \tag{59c}$$

From equations (47b) and (59), $\hat{\mathbf{t}}_s$ is given by

$$\hat{\mathbf{t}}_s = 2m_0 \hat{\mathbf{F}}_s^e \left( \frac{\partial \psi_{00}}{\partial I_1^e} \frac{\partial I_1^e}{\partial \hat{\mathbf{C}}^e} + \frac{\partial \psi_{00}}{\partial I_2^e} \frac{\partial I_2^e}{\partial \hat{\mathbf{C}}^e} + \frac{\partial \psi_{00}}{\partial I_3^e} \frac{\partial I_3^e}{\partial \hat{\mathbf{C}}^e} \right) \hat{\mathbf{F}}_s^{e\,\mathrm{T}}. \tag{60}$$

The yield surface $f$ is usually given by

$$f = (\hat{\mathbf{t}}_s : \mathbf{A}_1 : \hat{\mathbf{t}}_s)^{\frac{1}{2}}, \tag{61}$$

where $\mathbf{A}_1$ is a fourth-order tensor to describe the yielding behaviors. $f$ also can be considered as the effective stress $\bar{\sigma}$. $t_y$ is then given by

$$t_y = t_0 + H\Gamma_p, \tag{62}$$

where $t_0$ denotes initial yield stress and $H$ is a hardening slop related to $\frac{\partial t_{y0}}{\partial \Gamma_p}$ in equation (58). $\Gamma_p$ is the hardening variable and accumulated based on $D(\Gamma_p)$ in equation (58). The inelastic potential surface $Q_p$ can be given by

$$Q_p = (\hat{\mathbf{t}}_s : \mathbf{A}_2 : \hat{\mathbf{t}}_s)^{\frac{1}{2}}. \tag{63}$$

The coefficients in $\mathbf{A}_2$ should keep the convexity of $Q_p$ to satisfy the inequality condition of material dissipation in equation (47a).

In this example, an observer $A$ is defined at the origin of the frame $S$. In the two-dimensional Minkowski space, observer $A$ observes $\hat{\mathbf{x}}$, $\boldsymbol{\eta}$, $\hat{\mathbf{u}}$, $\boldsymbol{\Lambda}$, and $\mathbf{S}$ based on the equations in sections 2 and 3 as follows:

$$\hat{\mathbf{x}} = x^\alpha \hat{\mathbf{e}}_\alpha \ (a = 1-2), \text{ where } x^2 = ct, \tag{64a}$$

$$\boldsymbol{\eta} = \begin{bmatrix} 1 & 0 \\ 0 & -1 \end{bmatrix}, \tag{64b}$$

$$\hat{\mathbf{u}} = \begin{bmatrix} \frac{\beta}{\sqrt{1-\beta^2}} \\ \frac{1}{\sqrt{1-\beta^2}} \end{bmatrix}, \tag{64c}$$

$$\boldsymbol{\Lambda} = \begin{bmatrix} \frac{1}{\sqrt{1-\beta^2}} & -\frac{\beta}{\sqrt{1-\beta^2}} \\ -\frac{\beta}{\sqrt{1-\beta^2}} & \frac{1}{\sqrt{1-\beta^2}} \end{bmatrix}, \tag{64d}$$

$$\mathbf{S} = \begin{bmatrix} \frac{1}{1-\beta^2} & -\frac{\beta}{1-\beta^2} \\ \frac{\beta}{1-\beta^2} & -\frac{\beta^2}{1-\beta^2} \end{bmatrix}. \tag{64e}$$

In this example, $x^1$ and $x^2$ have unit-less values. With equations (7) and (64a), the contour on which the four-vectors have the same distance from the origin ($|\hat{\mathbf{x}}|^2 = l^2$) is calculated as

$$|\hat{\mathbf{x}}|^2 = \hat{\mathbf{x}} \cdot \hat{\mathbf{x}} = (x^1)^2 - (x^2)^2 = l^2. \tag{65}$$

Figure 2 plots the contour (dotted line in figure 2) and shows that the metric of the Minkowski space is different from that in the Euclidian space.

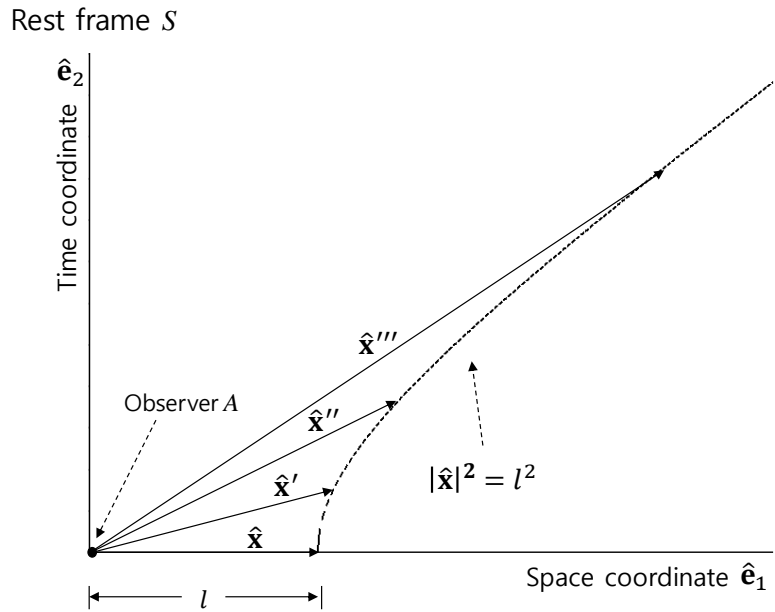

Figure 2. Contour of four vectors with same distances

A one-dimensional bar moving in the Minkowski space is set in figure 3. Figure 3 presents the world lines of two particles ($P_1(0)$ and $P_2(0)$ when $t = 0$) at the two ends of the bar observed by observer $A$. Length $L_0$ between $P_1(0)$ and $P_2(0)$ is the rest length of the bar observed by observer $A$ at $t = 0$. The position of a particle $x^1$, according to time is given by

$$x^1 = x^1(X, t), \tag{66a}$$

$$X = x^1 \text{ at } t=0. \tag{66b}$$

Thus, observer $A$ can observe $\hat{\mathbf{F}}_s$ based on equations (14), (16a), (64a) and (64e),

$$\hat{\mathbf{F}}_s = \mathbf{S} \begin{bmatrix} \frac{\partial x^1}{\partial X} \\ \frac{\partial x^2}{\partial X} \end{bmatrix} = \frac{1}{1-\beta^2} \begin{bmatrix} \frac{\partial x^1}{\partial X} - \beta \frac{\partial x^2}{\partial X} \\ \beta \frac{\partial x^1}{\partial X} - \beta^2 \frac{\partial x^2}{\partial X} \end{bmatrix}. \tag{67}$$

$\beta$ is defined by equation (10a). Thus, $\hat{\mathbf{C}}$, $\hat{\mathbf{B}}$, and $j$ are determined based on equations (17b), (19b), and

(67) by

$$\hat{C} = \left[\frac{1}{1-\beta^2}(\frac{\partial x^1}{\partial X} - \beta \frac{\partial x^2}{\partial X})^2\right], \tag{68a}$$

$$\hat{B} = \begin{bmatrix} \frac{1}{(1-\beta^2)^2}(\frac{\partial x^1}{\partial X} - \beta \frac{\partial x^2}{\partial X})^2 & \frac{\beta}{(1-\beta^2)^2}(\frac{\partial x^1}{\partial X} - \beta \frac{\partial x^2}{\partial X})^2 \\ \frac{\beta}{(1-\beta^2)^2}(\frac{\partial x^1}{\partial X} - \beta \frac{\partial x^2}{\partial X})^2 & \frac{\beta^2}{(1-\beta^2)^2}(\frac{\partial x^1}{\partial X} - \beta \frac{\partial x^2}{\partial X})^2 \end{bmatrix}, \tag{68b}$$

$$j = \frac{1}{\sqrt{1-\beta^2}}(\frac{\partial x^1}{\partial X} - \beta \frac{\partial x^2}{\partial X}). \tag{68a}$$

The invariants can be provided by

$$\hat{C}:I = \hat{B}:\eta I = \frac{1}{1-\beta^2}(\frac{\partial x^1}{\partial X} - \beta \frac{\partial x^2}{\partial X})^2, \tag{69a}$$

$$(\hat{C}I\hat{C}):I = (\hat{B}\eta\hat{B}):\eta I = (\frac{1}{1-\beta^2})^2(\frac{\partial x^1}{\partial X} - \beta \frac{\partial x^2}{\partial X})^4. \tag{69b}$$

As shown in figure 3, when $t = t_1$, observer $A$ simultaneously observes $P_1(t_1)$ and $P_2(t_2)$. Observer $A$ observes that the deformed length of the bar is $L$ at $t = t_1$. The change in length $L$ reflects not only the length change due to the deformation but also the relativistic effect, because the bar has a relative velocity between observer $A$; the distance between $P_1(t_1)$ and $P_2(t_2)$ is not the rest length. The relativity effect increases as the magnitude of the relative velocity approaches the speed of light based on $\beta$.

We define another frame $S'$ that is fixed to particle $P_1$, and the two move together. Another observer $A'$ is positioned at the origin of the frame $S'$. For observer $A'$, the time and space axes of frame $S'$ are geometrically independent; however, observer $A$ observes that frame $S'$ is distorted over space-time with $\hat{e}_1'$ and $\hat{e}_2'$ with respect to $\beta$, as shown in figure 3. Observer $A'$ instantaneously measures the rest length of the bar as $L'$ on the same timeline depending on $\hat{e}_1'$. $\hat{F}_s$, $\hat{C}$, and $\hat{B}$ representing the relationship between different observations by two observers. For observer $A$, deformed length $L$ can be obtained by

$$L = L_0 \hat{F}_{s(1,1)} = L_0\sqrt{\hat{B}_{(1,1)}} = \frac{L_0}{1-\beta^2}(\frac{\partial x^1}{\partial X} - \beta \frac{\partial x^2}{\partial X}). \tag{70}$$

$L'$ for observer $A'$ is obtained by

$$L' = L_0\sqrt{\hat{C}} = L_0\sqrt{\hat{B}:\eta I} = \frac{L_0}{\sqrt{1-\beta^2}}(\frac{\partial x^1}{\partial X} - \beta \frac{\partial x^2}{\partial X}). \tag{71}$$

The length in the time axis direction ($T$) is calculated as follows:

$$T = L_0 \hat{\mathrm{F}}_{s(2,1)} = L_0 \sqrt{\hat{\mathrm{B}}_{(2,2)}} = \frac{L_0}{1-\beta^2} \beta \left( \frac{\partial x^1}{\partial X} - \beta \frac{\partial x^2}{\partial X} \right). \tag{72}$$

Equations (70–72) show that $L$, $T$, and $L'$ satisfy the following relations:

$$L = \frac{1}{\sqrt{1-\beta^2}} L', \tag{73a}$$

$$T = \frac{\beta}{\sqrt{1-\beta^2}} L', \tag{73b}$$

$$L^2 - T^2 = L'^2. \tag{73c}$$

They physically mean that even though observer $A'$ instantaneously measures the deformed rest length of bar $L'$ in frame $S'$, observer $A$ experiences a time difference $\Delta t = \frac{T}{c}$ for $L'$. $\hat{\mathbf{B}}$ and $\hat{\mathbf{C}}$ present observations for observer $A$ and $A'$. However, based on equation (73c), both observers $A$ and $B$ can determine rest length $L'$, based on the invariance in equations (69a–69b). Figure 4 shows the relationship between $L$, $L'$, and $T$ observed by two observers ($A$ and $A'$), which is enlarged from figure 3.

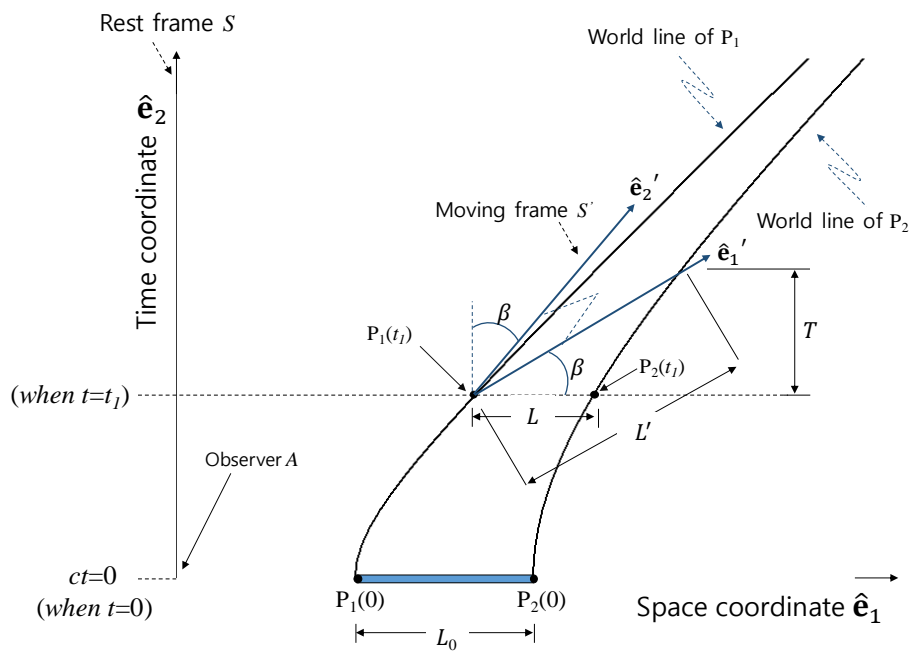

Figure 3. Physical interpretation of relativistic model in one bar example

https://doi.org/10.48550/arXiv.2305.12103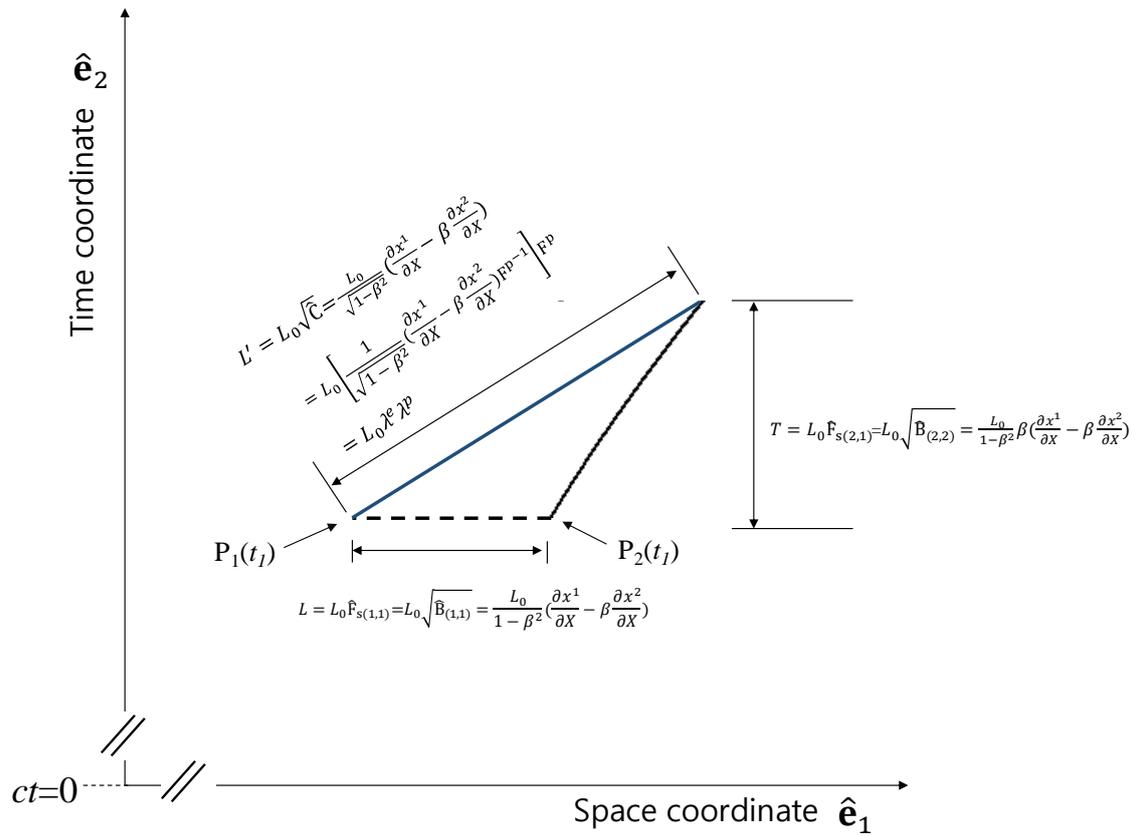

Figure 4. Relationship between different observations in one bar example

Next, the elastic and inelastic parts of the deformation gradient can be decomposed using

$$\hat{\mathbf{F}}^e = \frac{1}{(1-\beta^2)} \begin{bmatrix} \frac{\partial x^1}{\partial X} - \beta \frac{\partial x^2}{\partial X} \\ \beta \frac{\partial x^1}{\partial X} - \beta^2 \frac{\partial x^2}{\partial X} \end{bmatrix} F^{p-1}, \tag{74a}$$

$$F^p = \lambda^p. \tag{74b}$$

The elastic and inelastic parts of $\hat{\mathbf{B}}$ and $\hat{\mathbf{C}}$ are calculated based on equations (33) and (74),

$$\hat{\mathbf{B}}^e = \left[ \frac{(\frac{\partial x^1}{\partial X} - \beta \frac{\partial x^2}{\partial X})}{(1-\beta^2)} F^{p-1} \right]^2 \begin{bmatrix} 1 & \beta \\ \beta & \beta^2 \end{bmatrix}, \text{ and } \hat{\mathbf{B}}^p = \frac{\lambda^{p2}}{(1-\beta^2)} \begin{bmatrix} 1 & \beta \\ \beta & \beta^2 \end{bmatrix}, \tag{75a}$$

$$\hat{\mathbf{C}}^e = \frac{1}{(1-\beta^2)} \left[ (\frac{\partial x^1}{\partial X} - \beta \frac{\partial x^2}{\partial X}) F^{p-1} \right]^2, \text{ and } \hat{\mathbf{C}}^p = F^{p2}. \tag{75b}$$

Based on equations (31) and (64), $\hat{\mathbf{L}}_s^\eta$ and $\hat{\mathbf{D}}_s^\eta$ are calculated as follows:

$$\hat{\mathbf{L}}_s^\eta = (1-\beta^2)^{-\frac{5}{2}} \left( \frac{\partial \beta}{\partial x^1} + \beta \frac{\partial \beta}{\partial x^2} \right) \begin{bmatrix} 1 & -\beta \\ -\beta & \beta^2 \end{bmatrix}, \tag{76a}$$

$$\hat{\mathbf{D}}_s^\eta = (1-\beta^2)^{-\frac{5}{2}} \left(\frac{\partial \beta}{\partial x^1} + \beta \frac{\partial \beta}{\partial x^2}\right) \begin{bmatrix} 1 & -\beta \\ -\beta & \beta^2 \end{bmatrix}. \tag{76b}$$

To separate $\hat{\mathbf{L}}_s^\eta$ and $\hat{\mathbf{D}}_s^\eta$ into the inelastic and elastic parts, a relativistic constitutive equation should be specified. To calculate the stress, we define a simple isotropic free energy for this example as follows:

$$\psi_{00} = \frac{C_1}{2}(I_1^e - 1)^2, \text{ where } I_1^e = \frac{1}{(1-\beta^2)} \left[ \left(\frac{\partial x^1}{\partial X} - \beta \frac{\partial x^2}{\partial X}\right) F^{p-1} \right]^2. \tag{77}$$

Equations (60) and (76) yield the stress tensor,

$$\hat{\mathbf{t}}_s = 2m_0 C_1 \frac{\left(\frac{\partial x^1}{\partial X} - \beta \frac{\partial x^2}{\partial X}\right)^2}{(1-\beta^2)^3} F^{p-4} \left[\left(\frac{\partial x^1}{\partial X} - \beta \frac{\partial x^2}{\partial X}\right)^2 - \hat{\mathbf{C}}^p(1-\beta^2)\right] \begin{bmatrix} 1 & \beta \\ \beta & \beta^2 \end{bmatrix}, \tag{78a}$$

$$\hat{\mathbf{t}}_s^{(1,1)} = 2m_0 C_1 \frac{\left(\frac{\partial x^1}{\partial X} - \beta \frac{\partial x^2}{\partial X}\right)^2}{(1-\beta^2)^3} \hat{\mathbf{C}}^{p-2} \left[\left(\frac{\partial x^1}{\partial X} - \beta \frac{\partial x^2}{\partial X}\right)^2 - \hat{\mathbf{C}}^p(1-\beta^2)\right], \tag{78b}$$

$$\hat{\mathbf{t}}_s^{(2,2)} = \hat{\mathbf{t}}_s^{(1,1)} \beta^2, \tag{78c}$$

$$\hat{\mathbf{t}}_s^{(1,2)} = \hat{\mathbf{t}}_s^{(2,1)} = \beta \hat{\mathbf{t}}_s^{(1,1)}. \tag{78d}$$

In this example, the yield function is simply defined as

$$\bar{\sigma} = f = (\hat{\mathbf{t}}_{s(1,1)}^2 + \hat{\mathbf{t}}_{s(2,2)}^2 - \hat{\mathbf{t}}_{s(1,2)}^2 - \hat{\mathbf{t}}_{s(2,1)}^2)^{\frac{1}{2}}. \tag{79}$$

The inelastic potential is also defined as follows:

$$Q_p = (\hat{\mathbf{t}}_{s(1,1)}^2 + \hat{\mathbf{t}}_{s(2,2)}^2 - \hat{\mathbf{t}}_{s(1,2)}^2 - \hat{\mathbf{t}}_{s(2,1)}^2)^{\frac{1}{2}}. \tag{80}$$

When the yield surface $f$ satisfies the yielding condition in equation (54a), inelastic deformation occurs. Equations (32) and (52) provide $\frac{\partial f}{\partial \hat{\mathbf{t}}_s}$, $\overline{[\hat{\mathbf{D}}_s^\eta]}^p$, and $\overline{[\hat{\mathbf{L}}_s^\eta]}^p$,

$$\frac{\partial f}{\partial \hat{\mathbf{t}}_s} = \frac{\hat{\mathbf{t}}_s^{(1,1)}}{\bar{\sigma}} \begin{bmatrix} 1 & -\beta \\ -\beta & \beta^2 \end{bmatrix}, \tag{81a}$$

$$\overline{[\hat{\mathbf{D}}_s^\eta]}^p = \frac{\hat{\mathbf{t}}_s^{(1,1)}}{\bar{\sigma}} \begin{bmatrix} 1 & -\beta \\ -\beta & \beta^2 \end{bmatrix}, \tag{81b}$$

$$\overline{[\hat{\mathbf{L}}_s^\eta]}^p = \frac{\hat{\mathbf{t}}_s^{(1,1)}}{\bar{\sigma}} \begin{bmatrix} 1 & -\beta \\ -\beta & \beta^2 \end{bmatrix}. \tag{81c}$$

$[\hat{\mathbf{D}}_s^\eta]^p$ and $[\hat{\mathbf{L}}_s^\eta]^p$ become

$$[\hat{\mathbf{D}}_s^\eta]^p = \frac{\hat{\mathbf{t}}_s^{(1,1)}}{\bar{\sigma}} \begin{bmatrix} 1 & -\beta \\ -\beta & \beta^2 \end{bmatrix} D(\Gamma_p), \tag{82a}$$

$$[\hat{\mathbf{L}}_s^\eta]^p = \frac{\hat{t}_s^{(1,1)}}{\bar{\sigma}}\begin{bmatrix} 1 & -\beta \\ -\beta & \beta^2 \end{bmatrix} D(\Gamma_p). \tag{82b}$$

These tensors satisfy the inequality condition of the material by

$$\xi = \hat{\mathbf{t}}_s : [\widehat{\mathbf{D}}_s^\eta]^p = \frac{1}{\bar{\sigma}}\left[\hat{t}_s^{(1,1)}(1-\beta^2)D(\Gamma_p)\right]^2 \geq 0, \tag{83}$$

Equations (31) and (64) provide $[\widehat{\mathbf{D}}_s^\eta]^e$ and $[\hat{\mathbf{L}}_s^\eta]^e$,

$$[\widehat{\mathbf{D}}_s^\eta]^e = \widehat{\mathbf{D}}_s^\eta - [\widehat{\mathbf{D}}_s^\eta]^p = (1-\beta^2)^{-\frac{5}{2}}\left[\left(\frac{\partial \beta}{\partial x^1}+\beta\frac{\partial \beta}{\partial x^2}\right)-\frac{\hat{t}_{s(1,1)}}{\bar{\sigma}}(1-\beta^2)^{\frac{5}{2}}D(\Gamma_p)\right]\begin{bmatrix} 1 & -\beta \\ -\beta & \beta^2 \end{bmatrix}, \tag{84a}$$

$$[\hat{\mathbf{L}}_s^\eta]^e = \hat{\mathbf{L}}_s^\eta - [\hat{\mathbf{L}}_s^\eta]^p = (1-\beta^2)^{-\frac{5}{2}}\left[\left(\frac{\partial \beta}{\partial x^1}+\beta\frac{\partial \beta}{\partial x^2}\right)-\frac{\hat{t}_{s(1,1)}}{\bar{\sigma}}(1-\beta^2)^{\frac{5}{2}}D(\Gamma_p)\right]\begin{bmatrix} 1 & -\beta \\ -\beta & \beta^2 \end{bmatrix}. \tag{84b}$$

By substituting equations (74–84) into equation (59), $D(\Gamma_p)$ can be calculated as

$$D(\Gamma_p) = \frac{g_1}{g_2}, \text{ where} \tag{85a}$$

$$g_1 = 2\hat{t}_s^{(1,1)}(1-\beta^2)^{-\frac{3}{2}}(1+\beta^2)\left[\left(\frac{\partial \beta}{\partial x^1}+\beta\frac{\partial \beta}{\partial x^2}\right)-\frac{\hat{t}_{s(1,1)}}{\bar{\sigma}}(1-\beta^2)^{\frac{5}{2}}D(\Gamma_p)\right]+4m_0C_1(1-\beta^2)^{-\frac{7}{2}}\left(\frac{\partial \beta}{\partial x^1}+ \tag{85b}$$

$$\beta\frac{\partial \beta}{\partial x^2}\right)\left[\left(\frac{\partial x^1}{\partial X}-\beta\frac{\partial x^2}{\partial X}\right)F^{p-1}\right]^4,$$

$$g_2 = m_0\frac{\partial t_{y0}}{\partial \Gamma_p}+4m_0C_1(1-\beta^2)^{-2}\left[\left(\frac{\partial x^1}{\partial X}-\beta\frac{\partial x^2}{\partial X}\right)F^{p-1}\right]^4. \tag{85c}$$

Since equation (85) is a nonlinear function, it is effective to solve it numerically in most cases. The invariant derivative $D(\Gamma_p)$ is related to $\dot{\Gamma}_p$ as follows:

$$D(\Gamma_p) = \frac{1}{\sqrt{1-\beta^2}}\dot{\Gamma}_p, \text{ where} \tag{86a}$$

$$\dot{\Gamma}_p = \frac{2\hat{t}_s^{(1,1)}(1-\beta^2)^{-1}(1+\beta^2)\left[\left(\frac{\partial \beta}{\partial x^1}+\beta\frac{\partial \beta}{\partial x^2}\right)-\frac{\hat{t}_{s(1,1)}}{\bar{\sigma}}(1-\beta^2)^{\frac{5}{2}}D(\Gamma_p)\right]+4m_0C_1(1-\beta^2)^{-3}\left(\frac{\partial \beta}{\partial x^1}+\beta\frac{\partial \beta}{\partial x^2}\right)\left[\left(\frac{\partial x^1}{\partial X}-\beta\frac{\partial x^2}{\partial X}\right)F^{p-1}\right]^4}{m_0\frac{\partial t_{y0}}{\partial \Gamma_p}+4m_0C_1(1-\beta^2)^{-2}\left[\left(\frac{\partial x^1}{\partial X}-\beta\frac{\partial x^2}{\partial X}\right)F^{p-1}\right]^4}. \tag{86b}$$

$F^p$ can be calculated by

$$F^p = \exp\left[\int_{t_y}^{t_c}\dot{\Gamma}_p dt\right]. \tag{87}$$

$t_y$ and $t_c$ denote the start time of inelastic deformation and the current time, respectively. Finally, the objective elastic and inelastic deformation parts of the rest length $L'$, can be obtained, and the results are shown in figure 4.

$$L' = L_0\lambda^e\lambda^p, \tag{88a}$$

$$\lambda^e = \sqrt{\hat{\mathbf{C}}^e : \mathbf{I}} = \frac{1}{\sqrt{1-\beta^2}} \left( \frac{\partial x^1}{\partial X} - \beta \frac{\partial x^2}{\partial X} \right) \mathrm{F}^{\mathrm{p}-1}, \tag{88b}$$

$$\lambda^p = \sqrt{\hat{\mathbf{C}}^{\mathrm{p}} : \mathbf{I}} = \mathrm{F}^{\mathrm{p}}. \tag{88c}$$

### 4.2. Relationship with non-relativistic mechanics

This section discusses the relationship between the proposed relativistic constitutive model of inelastic deformation and the non-relativistic mechanics model for the same example. The link between the relativistic and nonrelativistic mechanical models is parameter $\beta = \frac{v}{c}$. If the effect on the relativity of the velocity is negligible ($v \ll c$), $\beta \approx 0$. In this case, $\hat{\mathbf{F}}_s$, $\hat{\mathbf{C}}$, $\hat{\mathbf{B}}$, and $j$ in section 4.1 are simplified as

$$\hat{\mathbf{F}}_s = \begin{bmatrix} \frac{\partial x^1}{\partial X} \\ 0 \end{bmatrix}, \tag{89a}$$

$$\hat{\mathbf{C}} = \left(\frac{\partial x^1}{\partial X}\right)^2, \tag{89b}$$

$$\hat{\mathbf{B}} = \begin{bmatrix} \left(\frac{\partial x^1}{\partial X}\right)^2 & 0 \\ 0 & 0 \end{bmatrix}, \tag{89c}$$

$$j = \frac{\partial x^1}{\partial X}. \tag{89d}$$

Setting $\beta = 0$ does not distort the time and space axes for either observer, and the variables related to the time axis affected by $\beta$ become zero. The surviving terms of $\hat{\mathbf{F}}_s$, $\hat{\mathbf{C}}$, and $\hat{\mathbf{B}}$ are the same as those of conventional non-relativistic tensors. In addition, the lengths of $L$, $T$, and $L'$ in section 4.1 become

$$L = L_0 \hat{\mathrm{F}}_{s(1,1)} = L_0 \sqrt{\hat{\mathrm{B}}_{(1,1)}} = L_0 \frac{\partial x^1}{\partial X}, \tag{90a}$$

$$T = L_0 \hat{\mathrm{F}}_{s(2,1)} = L_0 \sqrt{\hat{\mathrm{B}}_{(2,2)}} = 0, \tag{90b}$$

$$L_0' = L_0 \sqrt{\hat{\mathbf{C}}} = L_0 \sqrt{\hat{\mathbf{B}} : \mathbf{\eta}\mathbf{I}} = L_0 \frac{\partial x^1}{\partial X}. \tag{90c}$$

Because $L = L'$, the effects on the different velocities between frames $S$ and $S'$ disappear, similar to non-relativistic mechanics; time delay T becomes zero. The elastic and inelastic deformation parts are determined as follows:

$$L' = L_0 \lambda^e \lambda^p, \tag{91a}$$

$$\lambda^e = \sqrt{\widehat{\mathbf{C}^e \colon \mathbf{I}}} = \sqrt{\widehat{\mathbf{B}^e \colon \boldsymbol{\eta}\mathbf{I}}} = (\frac{\partial x^1}{\partial X})\mathrm{F}^{\mathrm{p}\,-1}, \tag{91b}$$

$$\lambda^p = \sqrt{\widehat{\mathbf{C}^p \colon \mathbf{I}}} = \sqrt{\widehat{\mathbf{B}^p \colon \boldsymbol{\eta}\mathbf{I}}} = \mathrm{F}^{\mathrm{p}}. \tag{91c}$$

These also reduce $[\widehat{\mathbf{D}}_s^{\boldsymbol{\eta}}]^{\mathrm{p}}$, $[\widehat{\mathbf{L}}_s^{\boldsymbol{\eta}}]^{\mathrm{p}}$, and $\hat{\mathbf{t}}_s$ to non-relativistic tensors as follows:

$$[\widehat{\mathbf{D}}_s^{\boldsymbol{\eta}}]^{\mathrm{p}} = \frac{\hat{\mathbf{t}}_s^{(1,1)}}{\bar{\sigma}}\begin{bmatrix}1 & 0\\ 0 & 0\end{bmatrix}\varGamma_p, \tag{92a}$$

$$[\widehat{\mathbf{L}}_s^{\boldsymbol{\eta}}]^{\mathrm{p}} = \frac{\hat{\mathbf{t}}_s^{(1,1)}}{\bar{\sigma}}\begin{bmatrix}1 & 0\\ 0 & 0\end{bmatrix}\dot{\varGamma}_p, \tag{92b}$$

$$\hat{\mathbf{t}}_s = \hat{\mathbf{t}}_s^{(1,1)}\begin{bmatrix}1 & 0\\ 0 & 0\end{bmatrix}. \tag{92c}$$

The results show that the proposed model is identical to non-relativistic model when the magnitude of the velocity is much smaller than the speed of light.

## 5. Conclusions

This paper presents a modeling structure of the relativistic constitutive equations of the inelastic deformation of a material moving in the Minkowski space. In this model, the extended four-dimensional Cauchy–Green deformation tensors, which are objective about the homogeneous Lorentz transformation in the Minkowski space-time, are used. Based on these four-dimensional tensors, constitutive equations of inelastic deformation are proposed to provide objective measures for different observers with different velocities over a wide speed range, even when the moving speed is close to the speed of light. The proposed constitutive modeling of inelastic deformation has the advantage of being able to use a previously introduced theory of inelasticity after modifications considering the relativity effect. Consequently, the proposed model satisfies both the traditional inelasticity theory and relativistic effect. In addition, the defined constitutive equations of relativistic inelastic deformation are shown to satisfy the second law of thermodynamics. The presented example shows how to apply the introduced theory by explaining the physical meanings of the formulations by comparing two observers in different frames. This shows that even though the two observers have different observations, they can have invariant measures for inelastic deformation. Finally, the proposed relativistic inelasticity model is proved to be identical to a nonrelativistic model when the speed of motion is much smaller than that of light.

# Conflict of interest

Authors declare there is no conflict of interest

# Acknowledgments

This work was supported from the National Research Foundation of Korea (NRF) grant (No. 2021R1C1C1007946) funded by the Korea government.